[1994/12/01]% LaTeX date must December 1994 or later
\documentclass[draft]{article}
\usepackage{amsmath}
\usepackage{amsfonts}
\usepackage{amssymb}
\usepackage{euscript}
\input epsf
\usepackage{amsthm}
\newtheorem{thm}{Theorem}[section]

\newtheorem{lem}[thm]{Lemma}

\theoremstyle{definition}
\newtheorem{defn}{Definition}[section]

\theoremstyle{remark}
\newtheorem{rem}{Remark}[section]

\begin{document}
\begin{center}

 %PLEASE INSERT THE TITLE OF YOUR CONTRIBUTION HERE
{\large\bf
A COMPARISON BETWEEN \\
DIFFERENT CONCEPTS  OF  \\
ALMOST ORTHOGONAL POLYNOMIALS
 }
\bigskip

 % PLEASE INSERT THE NAMES AND ADDRESSES OF THE AUTHORS HERE
\centerline{\bf Predrag Rajkovi\'c${}^1$ and Sladjana Marinkovi\'c${}^2$ }
% and Sladjana Marinkovi\'c${}^1$
\end{center}

\medskip

\centerline{${}^1$Department of Mathematics, Faculty of Mechanical
Engineering,}

\centerline{University of Ni\v s, A. Medvedeva 14, 18 000 Ni\v s,
Serbia}

\centerline{\tt pedja.rajk@masfak.ni.ac.rs}

\smallskip

\centerline{${}^2$Department of Mathematics, Faculty of
Electronic Engineering,}

\centerline{\tt sladjana@elfak.ni.ac.rs}

\bigskip

% PLEASE INSERT THE ABSTRACT HERE
\noindent {\bf Abstract.} In this paper, we will discuss the notion of almost orthogonality in a functional sequence.Especially, we will define a few sequences of almost orthogonal polynomials which can be used successfully for modeling of electronic systems which generate orthonormal basis. We will include quasi-orthogonality and examine its influence on the behavior of these sequences.

\medskip

%42C05 Orthogonal functions and polynomials, general theory [See
%also 33C45, 33C50, 33D45]

%33C45 Orthogonal polynomials and functions of hypergeometric type
%(Jacobi, Laguerre, Hermite, Askey scheme, etc.) [See also 42C05
%for general orthogonal polynomials and functions]

%05E35 Orthogonal polynomials [See also 33C45, 33C50, 33D45]

\noindent  {\it Keywords.} Operator, Functional, Function, Polynomial, Orthogonality, Quasi orthogonality, Almost orthogonality.

\smallskip

\noindent {\it Mathematics Subject Classification.} Primary 42C05,
Secondary 33C45.
\bigskip

% Section 1
\section{Introduction}

The first usage of the notion of {\it almost orthogonality} for operators is annotated in the M. Cotlar's paper \cite{Cotlar}. Let $E$ and $F$ be the Hilbert spaces with their scalar products and norms. For a linear operator $S:\ E\to F$, the operator $S^\star:\ F\to E$ is {his adjoined operator} if it is satisfied
\begin{equation}\label{Adjoint}
(Su,v)_F=(u,S^\star v)_E \qquad (\forall u\in E,\ \forall v\in F).
\end{equation}
The operator norm is
\begin{equation}\label{NormOp}
\|S\|=\sup_{\|u\|_E=1} \| S(u)\|_F, \qquad
\|S^\star\|=\sup_{\|v \|_F=1} \| S^\star(v)\|_E.
\end{equation}

\begin{defn} (Almost orthogonal operators). We will call a family of continuous operators
%\begin{equation}\label{Cotlar1}
$T_i:\ E\to F\quad (i\in\mathbb Z),$
%\end{equation}
{\it almost orthogonal} if they satisfy the following conditions:
\begin{equation}\label{Cotlar2}
\|T^\star_iT_j\|\le a_{i,j}, \quad \|T_iT^\star_j\|\le b_{i,j}, \qquad (i,j\in\mathbb Z),
\end{equation}
where $a_{i,j}$  and $b_{i,j}$ are non-negative symmetric functions on $\mathbb Z\times \mathbb Z$ which satisfy
\begin{equation}\label{Cotlar3}
\|a\|^\mu_{\infty,\mu}=\sup_{i\in\mathbb Z} \sum_{j\in\mathbb Z} a^{\mu}_{i,j}<\infty \quad
\|b\|^\nu_{\infty,\nu}=\sup_{i\in\mathbb Z} \sum_{j\in\mathbb Z} b^{\nu}_{i,j} <\infty, \end{equation}
where
%\begin{equation}\label{Cotlar4}
$0\le \mu,\nu \le 1,\quad  \mu+\nu=1.$
%\end{equation}
\end{defn}

\begin{lem}(Cotlar-Stein Lemma). Let $\{T_i\}_{i\in\mathbb Z}$ be a family of almost orthogonal operators. Then the formal sum $\sum_{i} T_i$ converges in the strong operator topology to a continuous linear operator $T : E \to F$,
which is bounded by
\begin{equation}\label{Cotlar5}
\|T\|\le \sqrt{\|a\|^\mu_{\infty,\mu} \ \|b\|^\nu_{\infty,\nu}}\ .
\end{equation}
\end{lem}

The concept of quasi-orthogonality was introduced in 1923. by M. Riesz \cite{Riesz} who considered the moment problem.
%Riesz, M. { Sur le problµeme des moments, III, Ark. Mat. Astr. Fys., 17(16) (1923),
%1--52.
It also appeared in Fejer's research of quadratures \cite{Fejer} in 1933.
%Fejer, L. { Mechanischen quadraturen mit positiven Cotesschen Zahlen, Math.
%Zeitschrift, 37 (1933), 287--309.
Later, a various aspects of this theory were considered by other mathematicians (T.S. Chihara \cite{ChiharaQuasi}, D.J. Dickinson \cite{Dickinson}, F. Marcelan,\ldots).

\begin{defn} (Quasi orthogonal functions).
We say that a  functional sequence $\{Q_n(x)\}$ is {\it quasi-orthogonal} of order $\rho$ $(\rho\in\mathbb N_0)$ with respect to the functional $U$ if
\begin{equation}\label{QuasiOPG}
U[Q_mQ_n]=0 \qquad (m,n\in\mathbb N_0:|m-n|>\rho).
\end{equation}
In the special case $\rho=0$, it becomes the regular orthogonality.
\end{defn}

In our paper \cite{Dankovic}, we have introduced the next concept.
\begin{defn} (Almost orthogonality by an error matrix)
Let  $\mathcal E=[\varepsilon_{i,j}]$ be a matrix whose elements are very small positive real numbers. If it exists, the sequence of the functions
$\{P^{(\varepsilon)}_n(x)\}$ which satisfies the relation
\begin{equation}\label{Al4}
\mathcal L\Bigl[P^{(\varepsilon)}_n \cdot
P^{(\varepsilon)}_i\Bigr]=\varepsilon_{n,i} \quad (i=0,1,\ldots,n-1;\ n\in
\mathbb N)
\end{equation}
 will be called {\it almost orthogonal} with respect to
$\mathcal L$ and the error matrix $\mathcal E$.
\end{defn}

\section{Almost orthogonality by shifted zeros}

Let $\lambda(x)$ be a positive Borel measure on an interval
$(a,b)\subset \mathbb R$ with infinite support and such that all
moments
\begin{equation}\label{Al1}
\lambda_n=\mathcal L[x^n] = \int_a^b x^n d\lambda (x)
\end{equation}
exist. In this manner, we define linear functional $\mathcal L$ in
the linear space of real polynomials $\mathcal P$. Also, we can
introduce an inner product as follows (see \cite{Szego}):
\begin{equation}\label{Al2}
(f,\ g) = \mathcal L\bigl[f\cdot g\bigr] \qquad (f,g\in\mathcal
P),
\end{equation}
which is positive-definite because of the property
$\|f\|^2=(f,f)\ge 0$. Hence it follows  that monic polynomials
$\{ P_n(x)\}$ orthogonal with respect to this inner product exist
and they satisfy the three-term recurrence relation
\begin{equation}\label{RR3}
 P_{k+1}(x) = (x - \alpha_k) P_{k}(x) - \beta_k  P_{k-1}(x) \quad
(k\ge 0),  \qquad  P_{-1} \equiv 0,\  P_{0} \equiv 1.
\end{equation}
The zeros of these polynomials are all contained in the interval $(a,b)=\textrm{supp}\lambda (x)$ and they interlace each other. If we denote them by $\{x_{n,k}\}$, we can write
\begin{equation}\label{OP1}
P_{n}(x) = \prod_{k=1}^n (x-x_{n,k}).
\end{equation}

Let us denote by
\begin{equation}\label{OP3}
\tilde P_{n}(x) = \sigma_n P_{n}(x),\quad \textrm{where}\quad \sigma_n=\frac{1}{\| P_n\|} .
\end{equation}
Obviously, $\{\tilde P_{n}(x)\}$ is the corresponding orthonormal polynomial
sequence.
\begin{equation}\label{OrtoNormRR3}
\sqrt{\beta_{n+1}} \tilde P_{n+1}(x) = (x - \alpha_n) \tilde P_{n}(x) - \sqrt{\beta_n}  \tilde P_{n-1}(x) \qquad (n\ge 0),
\end{equation}
\begin{equation}\label{OrtoNormRR3_1}
 \tilde P_{-1} \equiv 0,\quad  \tilde P_{0} \equiv \frac{1}{\sqrt{\beta_0}}.
\end{equation}

The next lemma, proven in  \cite{Brent}, will be very useful

\begin{lem}\label{PrMinors}  All leading principal minors of the matrix
\begin{equation}\label{PrinMinors}
A=\bmatrix
P_{n-1}(x_{n,1})& P_{n-1}(x_{n,2}) & \cdots & P_{n-1}(x_{n,n})\\
P_{n-2}(x_{n,1}) & P_{n-2}(x_{n,2}) &        & P_{n-2}(x_{n,n}) \\
\vdots \\
P_{0}(x_{n,1}) & P_{0}(x_{n,1}) & \cdots & P_{0}(x_{n,1})\\
\endbmatrix
\end{equation}
are nonsingular.
\end{lem}

Let us remind on notation
\begin{equation}\label{O1}
\alpha(x)=\mathcal O(\varepsilon^\beta)\quad \Leftrightarrow\quad \lim_{\varepsilon\to 0} \frac{\alpha(x)}{\varepsilon^\beta}=c \quad(0<c<\infty).
\end{equation}

The next lemma is slightly generalization of similar one from \cite{Wilkinson}.
\begin{lem} \label{Lemma1.2}Let $z_r$ be an isolated zero of a polynomial $f(z)$ and $g(z)$ a continuous function in $z_r$. Then the function
\begin{equation}\label{O2}
T(z)=f(z)+\varepsilon g(z) \quad(0<\varepsilon\ll 1)
\end{equation}
has a zero $z_r(\varepsilon)$ such that
\begin{equation}\label{O3}
z_r(\varepsilon) = z_r -\varepsilon \frac{g(z_r)}{f^\prime(z_r)}+\mathcal O(\varepsilon^2).
\end{equation}
\end{lem}
\noindent{\it Proof.} Under assumptions, we have $f(z_r)=0$, $f^\prime(z_r)=\kappa\ne0$ and
$T(z_r(\varepsilon))=0$.
According to mean valued theorem, we can write
\begin{equation*}%\label{O4}
\frac{f(z_r(\varepsilon))-f(z_r)}{z_r(\varepsilon) - z_r} = f^{\prime}(\eta_r(\varepsilon)) \quad (\eta_r(\varepsilon)\in (\min\{z_r(\varepsilon), z_r\}, \ \max\{z_r(\varepsilon), z_r\})),
\end{equation*}
Including it into (\ref{O2}), we have
\begin{equation*}%\label{O4}
T(z_r(\varepsilon))= f^{\prime}(\eta_r(\varepsilon))(z_r(\varepsilon)- z_r)+\varepsilon \ g(z_r(\varepsilon))=0,
\end{equation*}
wherefrom
\begin{equation}\label{O5}
z_r(\varepsilon) - z_r = -\varepsilon \ \frac{g(z_r(\varepsilon))}{f^{\prime}(\eta_r(\varepsilon))}.
\end{equation}
Since $f^{\prime}$ and $g$ are continuous functions in the point $z_r$, we can write
\begin{equation*}%\label{O7}
\varphi(\varepsilon) = \frac{f^{\prime}(\eta_r(\varepsilon))}{g(z_r(\varepsilon))}
=\frac{f^{\prime}(z_r)+k_1\varepsilon}{ g(z_r)+k_2\varepsilon}.
\end{equation*}
By using Taylor series of the function $\varphi(\varepsilon)$, we obtain
\begin{equation*}%\label{O7}
\varphi(\varepsilon) =
\frac{f^{\prime}(z_r)}{ g(z_r)}+\varepsilon\frac{k_1f^{\prime}(z_r)-k_2g(z_r)}{ [g(z_r)]^2}+\mathcal O(\varepsilon^2).
\end{equation*}
Hence we finish the proof of the formula (\ref{O3}). $\square$

For the next two theorems we find inspiration in R. Brent's paper \cite{Brent}. There, discussion about almost orthogonality was motivated by iterative methods for zero-finding, but we find that echo of this paper could be large in the theory of orthogonality itself. Our purpose is to improve conclusions in that way.

Let
\begin{equation}\label{Gamma1}
0<\varepsilon\ll 1,\quad s\in\{1,\ldots, n-1\},\quad |\gamma_{n,k} - x_{n,k} |<\varepsilon \qquad (k=1,\ldots,s),
\end{equation}
where $x_{n,k}$ are the zeros of $P_n(x)$ given by (\ref{OP1}).

\begin{thm} \label{Lemma1.3}
Under the condition (\ref{Gamma1}), the polynomial
\begin{equation}\label{q1}
Q_n(x)=\sigma_n \prod_{i=1}^s (x-\gamma_{n,i})\prod_{i=s+1}^n
(x-x_{n,i}),
\end{equation}
is almost orthogonal with respect to $\{\tilde
P_k(x)\}_{k=0}^n$, i.e.
\begin{equation}\label{Almost4}
f_i = \mathcal L\bigl[\tilde P_{i} \ Q_n\bigr] =\left\{
\begin{array}{cc}
\varepsilon\ \omega_i \ ,&0\le i\le n-1\ ,\\  \\
1\ ,& i=n\ ,
\end{array}
\right. \quad (\omega_i\in\mathbb R,\ 1\le i \le n).
\end{equation}
\end{thm}
\noindent{\it Proof.}
Let us denote by
$$
R_{n;k_1,k_2,\ldots,k_\ell}^{(\ell)}(x)= \frac{\tilde P_n(x)}{\displaystyle{
\prod_{i=1}^\ell (x-x_{n,k_i})}}\quad
(1\le k_1< \ldots < k_\ell\le n,\ 1\le \ell\le n,\ n\in\mathbb N).
$$
According to (\ref{Gamma1}), we can write
\begin{equation}\label{Gamma2}
\gamma_{n,k} = x_{n,k} +\varepsilon_{n,k}, \quad
\textrm{where}\quad |\varepsilon_{n,k}|<\varepsilon \quad
(k=1,\ldots,s).
\end{equation}
Then
\begin{equation}\label{Almost30}
Q_n(x) = \tilde P_n(x) + \sum_{m=1}^s (-1)^m\sum_{1\le i_1<\cdots <i_m \le s}
\prod_{k=1}^m\varepsilon_{n,i_k} \  R_{n;i_1,i_2,\ldots,i_m}^{(m)}(x) .
\end{equation}
Hence
\begin{equation}\label{Almost3}
Q_n(x)= \tilde P_n(x)+R_{n-1}(x)\mathcal O(\varepsilon)\qquad
(R_{n-1}\in\mathcal P),
\end{equation}
wherefrom the conclusion follows. $\square$.

\medskip

Especially, let be
\begin{equation}\label{Almost311}
R_{n,k}(x)\equiv R_{n,k}^{(1)}(x)=\frac{\tilde
P_n(x)}{x-x_{n,k}},\quad \tau_{i,n,k}=\mathcal L\bigl[\tilde P_{i}
\ R_{n,k}\bigr] \quad (1\le k \le n).
\end{equation}

Because of orthogonality, we can write
\begin{equation*}%\label{Almost312}
0 = \mathcal L\bigl[\tilde P_{i} \ \tilde P_{n}\bigr]=\mathcal L\bigl[\tilde P_{i}
\ (x-x_{n,k})R_{n,k}\bigr]=\mathcal L\bigl[x\tilde P_{i} \
R_{n,k}\bigr]-x_{n,k}\mathcal L\bigl[\tilde P_{i} \ R_{n,k}\bigr].
\end{equation*}
From three-term recurrence relation (\ref{OrtoNormRR3}), we have
\begin{equation}\label{Almost312}
x\tilde P_{i}(x) = \sqrt{\beta_{i+1}}\tilde P_{i+1}(x)
+ \alpha_i \tilde P_{i} (x)
+ \sqrt{\beta_{i}} \tilde P_{i-1}(x).
\end{equation}
Hence
\begin{equation}\label{OrtoNormRR31}
\sqrt{\beta_{i+1}}\ \tau_{i+1,n,k} = (x_{n,k} - \alpha_i)\ \tau_{i,n,k} - \sqrt{\beta_i}
 \ \tau_{i-1,n,k} \quad (0\le i<n; 1\le k \le n).
\end{equation}

\begin{lem} \label{tau}
Let
\begin{equation}\label{Cond31}
h = \min_{0\le i \le n}\sqrt{\beta_{i}},\quad
R = \max_{0\le i \le n}\sqrt{\beta_{i}},\quad
C = \max_{{0\le i \le n}\atop{1\le k \le n}} |x_{n,k}-\alpha_i|.
\end{equation}
Then
\begin{equation}\label{Cond32}
|\tau_{i,n,k} | \le |\tau_{0,n,k} |\Bigl(\frac{C}{h}\Bigr)^i\sum_{j=0}^{[i/2]} {{i-j}\choose {j}} \Bigl(\frac{Rh}{C^2}\Bigr)^j.
\end{equation}
\end{lem}
\noindent{\it Proof.} By mathematical induction.$\square$

\medskip

By using the form (\ref{Almost30}) of the polynomial $Q_n(x)$, we can write
$$
\aligned
&f_i = \mathcal L\bigl[\tilde P_{i} \ Q_n\bigr]
=\mathcal L\bigl[\tilde P_{i} \ \tilde P_n(x)\bigr]
-\sum_{k=1}^s \varepsilon_{n,k} \mathcal L\bigl[\tilde P_{i} \ R_{n,k}\bigr]\\
&+\sum_{{k_1,k_2=1}\atop{k_1<k_2}}^s \varepsilon_{n,k_1}\varepsilon_{n,k_2} \mathcal L\bigl[\tilde P_{i} \ R_{n;k_1,k_2}^{(2)}\bigr]
 \ +\ \cdots+(-1)^s \mathcal L\bigl[\tilde P_{i} \ R_{n;1,2,\ldots,s}^{(s)}\bigr] \prod_{i=1}^s\varepsilon_{n,i} \ .
\endaligned
$$
Hence
$$
f_i=-\sum_{k=1}^s \varepsilon_{n,k}\ \tau_{i,n,k}+\mathcal O(\varepsilon^2)\quad (1\le i\le n-1),\qquad f_n=1.
$$
According to (\ref{Gamma2}) and Lemma \ref{tau}, the following estimate is valid:
$$
|f_i|\le \varepsilon s |\tau_{i,n,k} |+\mathcal O(\varepsilon^2)\ .
$$
So, we can say that
$$
\omega_i\le s  |\tau_{0,n,k} | \Bigl(\frac{C}{h}\Bigr)^i\sum_{j=0}^{[i/2]}
{{i-j}\choose {j}} \Bigl(\frac{Rh}{C^2}\Bigr)^j+\mathcal
O(\varepsilon^2)\qquad (i=0,1,\ldots,n-1)\ .
$$

Notice that
\begin{equation*}%\label{Almost4}
Q_n(x)=\sum_{i=0}^n f_i \tilde P_i(x).
\end{equation*}
\begin{thm}
Under the condition (\ref{Gamma1}), the real numbers $\gamma_{n,s+1},\ldots, \gamma_{n,n}$ exist such that
\begin{equation*}%\label{Almost1}
\gamma_{n,k} = x_{n,k} + \mathcal O(\varepsilon)\qquad (k=s+1,\ldots,n),
\end{equation*}
and the polynomial
\begin{equation*}%\label{Almost1}
P^{(\varepsilon)}_n(x)=\sigma_n \ \prod_{k=1}^n (x-\gamma_{n,k})
\end{equation*}
is quasi almost orthogonal with respect to $\{P_k(x)\}_{k=0}^n$, i.e.
\begin{equation*}%\label{Almost2}
\mathcal L\bigl[\tilde P_{k}  P^{(\varepsilon)}_n\bigr] =  \left\{
\begin{array}{lc}
0\ ,                      &0\le k\le n-s-1\ ,\\  \\
\mathcal O(\varepsilon)\ ,& n-s\le k\le n-1\ ,\\  \\
1\ ,                      &k= n\ .
\end{array}
\right.
\end{equation*}
\end{thm}
\noindent{\it Proof.} Using the same notation like in the previous lemmas, we can define
\begin{equation}\label{q20}
T_n(x)=Q_n(x)+\varepsilon\left\{-\sum_{i=0}^{n-s-1} \omega_i \tilde P_i(x)+\sum_{i=n-s}^{n-1} \mu_i \tilde P_i(x)\right\},
\end{equation}
where constants $\mu_i$ \ $(i=n-s,\ldots,n-1)$ will be determined.
Also, it can be written in the form
\begin{equation}\label{q2}
T_n(x) = \tilde P_n(x) + \varepsilon\sum_{i=n-s}^{n-1} (\omega_i+\mu_i) \tilde P_i(x).
\end{equation}
Then we find
\begin{equation}\label{gj1}
g_j=\mathcal L\bigl[\tilde P_{j} \ T_n\bigr]=\left\{
\begin{array}{cc}
0\ ,&0\le j\le n-s-1\ ,\\  \\
\varepsilon(\omega_j+\mu_j)\ ,&n-s\le j\le n-1\ ,\\ \\
1\ ,& j=n\ .
\end{array}
\right.
\end{equation}
If $\{t_{n,k}\}$ are the zeros of $T_n(x)$, we can write
\begin{equation}\label{q21}
T_n(x) = \sigma_n\prod_{k=1}^n (x-t_{n,k}).
\end{equation}
By applying Lemma \ref{Lemma1.2} onto (\ref{q2}), for $k=1,\ldots, s$, we can write
\begin{equation}\label{Gamma''}
t_{n,k} = \gamma_{n,k} + \varepsilon\left\{\sum_{i=0}^{n-s-1} \omega_i\ \frac{
\tilde P_i(\gamma_{n,k})}{Q'_n(\gamma_{n,k})} - \sum_{i=n-s}^{n-1} \mu_i
\frac{\tilde  P_i(\gamma_{n,k})}{Q'_n(\gamma_{n,k})}\right\}+\mathcal
O(\varepsilon^2) .
\end{equation}
It can be written in the matrix form
\begin{equation}\label{q3}
A_s(\varepsilon) \vec \mu= \vec b(\varepsilon),
\end{equation}
where
\begin{equation}\label{PMinors}
A_s(\varepsilon)=\bmatrix
\tilde P_{n-s}(\gamma_{n,1}) & \cdots &\tilde P_{n-1}(\gamma_{n,1}) \\
\vdots \\
 \tilde P_{n-s}(\gamma_{n,s})&  & \tilde P_{n-1}(\gamma_{n,s})
\endbmatrix,\
\vec \mu=\bmatrix
\mu_{n-s}\\
\vdots \\
\mu_{n-1}
\endbmatrix,
\
\vec b(\varepsilon)=\bmatrix
b_1(\varepsilon)\\
\vdots \\
b_s(\varepsilon)
\endbmatrix,
\end{equation}
with
\begin{equation}\label{b}
b_k(\varepsilon)= Q^\prime_n(\gamma_{n,k})\frac{\gamma_{n,k}-t_{n,k}}{\varepsilon}+\sum_{j=0}^{n-s-1} \omega_j\  \tilde P_j(\gamma_{n,k})
+\mathcal O(\varepsilon)\qquad (k=1,\ldots, s).
\end{equation}
Let us consider the system
\begin{equation}\label{system}
A_s(\varepsilon) \vec \mu= \vec b^\prime(\varepsilon),\quad \text{\rm where}
\quad
b^\prime_k(\varepsilon)=\sum_{j=0}^{n-s-1} \omega_j \tilde P_j(\gamma_{n,k})+\mathcal O(\varepsilon).
\end{equation}
According to Lemma~\ref{PrMinors}, all leading principal minors of the matrix $A$, defined by (\ref{PrinMinors}), are nonsingular. Hence, for sufficiently small $\varepsilon$, the matrix $A_s(\varepsilon)$ is
nonsingular too. Therefore exists the solution
$$
\vec \mu=A_s^{-1}(\varepsilon)\ \vec b^\prime(\varepsilon).
$$
of the system (\ref{system}). In that case, it is valid $t_{n,k}=\gamma_{n,k}$ $(k=1,\ldots , s)$.

Taking
$\gamma_{n,k}=t_{n,k}$ $(k=s+1,\ldots ,n)$, we have
\begin{equation}\label{q22}
T_n(x) = \sigma_n\prod_{k=1}^n (x-\gamma_{n,k}),
\end{equation}
and
$$
\gamma_{n,k}=x_{n,k}+\mathcal O(\varepsilon)\qquad (k=1,2,\ldots n).
$$
Choosing $P^{(\varepsilon)}_n(x)=T_n(x)$, we prove its existence. $\square$

\begin{rem}
Because of its quasi orthogonality, the sequence  $\{P^{(\varepsilon)}_n(x)\}$ satisfies $(s+2)$-term recurrence relation of the form
\begin{equation}\label{AlmostRecsplus2}
x\ P^{(\varepsilon)}_n(x) =\sum_{k=n-s}^{n+1} d_{n,k}\ P^{(\varepsilon)}_k(x).
\end{equation}
\end{rem}

\begin{rem}
Writing $P^{(\varepsilon)}_n(x)$ in the form
\begin{equation*}
P^{(\varepsilon)}_n(x)= \sigma_n\prod_{k=1}^s (x-\gamma_{n,k})\ V_{n-s}(x),\quad \text{\rm where}\quad V_{n-s}(x)=x^{n-s}+\sum_{i=0}^{n-s-1} v_{n-s,i}x^i,
\end{equation*}
we can evaluate numerically $v_{n-s,i}$ $\ (0\le i\le n-s-1)$ from linear algebraic system obtained from the fact
$$
\mathcal L\bigl[P_{j} \ P^{(\varepsilon)}_n\bigr]=0 \quad (0\le j\le n-s-1).
$$
\end{rem}

\subsection{Examples}

In the examples we will take upper limit $\varepsilon$ and choose $\varepsilon_{n,k}$ by the function \text{\tt Random} from package {\it Mathematica} in the interval $(-\varepsilon, \varepsilon)$. We will repeat the whole procedure $20$ times.
Let us consider $P_4(x)$ and $P^{(\varepsilon)}_4(x)$ provided by $s=2$.
In the tables, the notation $a(n)$ means $a\cdot 10^{n}$.
In the first column is $\varepsilon$. In the second column is the maximal distance between the zeros of orthogonal and almost orthogonal polynomials of the same degree. In the third column, it is given maximal absolute value of inner products $\mathcal L\bigl[P_{j} \ P^{(\varepsilon)}_n\bigr]$ $(j=0,1,\ldots, n-1)$, some kind of almost orthogonality between the members.

\smallskip

\textbf{Example 1}. Let us consider Legendre polynomials $\{P_n(x)\}$  which are orthonormal with respect to the functional
$$
\mathcal L[f\cdot g]=\int_{-1}^{1}f(x)g(x)\,dx\ .
$$

\begin{center}
\begin{tabular}[b]{|c|c|}\hline
{\rm Weak}  &  {\rm orthogonality}\rule[-1.5ex]{0pt}{4ex}\\ \hline
$\varepsilon$ & {\rm inner products}\rule[-1.5ex]{0pt}{4ex}\\ \hline
0.1(-1) & 0.377056 (-1) \rule[-1.5ex]{0pt}{6ex} \\
0.1(-2) & 0.446765 (-2) \rule[-1.5ex]{0pt}{3ex} \\
0.1(-3) & 0.158524 (-2) \rule[-1.5ex]{0pt}{3ex} \\
0.1(-4) & 0.191905 (-3) \rule[-1.5ex]{0pt}{3ex} \\
\hline
\end{tabular}
\qquad
\begin{tabular}[b]{|c|c|c|}\hline
{\rm Quasi} &  {\rm almost} &  {\rm orthogonality}\rule[-1.5ex]{0pt}{4ex}\\ \hline
$\varepsilon$ & {\rm zero distance}& {\rm inner products}\rule[-1.5ex]{0pt}{4ex}\\ \hline
0.1(-1) & 0.500665(-1) & 0.186743(0) \rule[-1.5ex]{0pt}{6ex} \\
0.1(-2) & 0.584278(-2) & 0.181221(-1) \rule[-1.5ex]{0pt}{3ex} \\
0.1(-3) & 0.493928(-3) & 0.158524(-2) \rule[-1.5ex]{0pt}{3ex} \\
0.1(-4) & 0.602947(-4) & 0.191905(-3) \rule[-1.5ex]{0pt}{3ex} \\
\hline
\end{tabular}
\end{center}
\medskip

\textbf{Example 2}. The Laguerre polynomials $\{L_n(x)\}$ are orthonormal with respect to the functional
$$
\mathcal L[f\cdot g]=\int_0^{+\infty}f(x)g(x){\rm e}^{-x}\,dx\ .
$$

\begin{center}
\begin{tabular}[b]{|c|c|}\hline
{\rm Weak} &  {\rm orthogonality}\rule[-1.5ex]{0pt}{4ex}\\ \hline
$\varepsilon$ &  {\rm inner products}\rule[-1.5ex]{0pt}{4ex}\\ \hline
0.1(-1) &  0.166269(-1) \rule[-1.5ex]{0pt}{6ex} \\
0.1(-2) &  0.145173(-2) \rule[-1.5ex]{0pt}{3ex} \\
0.1(-3) &  0.174978(-3) \rule[-1.5ex]{0pt}{3ex} \\
0.1(-4) &  0.177084(-4) \rule[-1.5ex]{0pt}{3ex} \\
\hline
\end{tabular}
\qquad
\begin{tabular}[b]{|c|c|c|}\hline
{\rm Quasi} &  {\rm almost} &  {\rm orthogonality}\rule[-1.5ex]{0pt}{4ex}\\ \hline
$\varepsilon$ & {\rm zero distance}& {\rm inner products}\rule[-1.5ex]{0pt}{4ex}\\ \hline
0.1(-1) & 0.231992(1) & 0.647862(0) \rule[-1.5ex]{0pt}{6ex} \\
0.1(-2) & 0.167538(0) & 0.499177(-1) \rule[-1.5ex]{0pt}{3ex} \\
0.1(-3) & 0.178939(-1) & 0.530434(-2) \rule[-1.5ex]{0pt}{3ex} \\
0.1(-4) & 0.167632(-2) & 0.497250(-3) \rule[-1.5ex]{0pt}{3ex} \\
\hline
\end{tabular}
\end{center}

We can notice from the tables that insisting  to have quasi orthogonality included, has the consequence weakness of almost orthogonality of the members with high degrees, i.e. increasing of the values of inner products $\mathcal L\bigl[P_{j} \ P^{(\varepsilon)}_n\bigr]$ for high $j$'s.

\smallskip

\noindent\noindent {\bf Acknowledgements.} The research of the
authors was supported by Ministry of Science of Republic Serbia
through the Project 144023.

\end{document}